\documentclass{article}
\usepackage[utf8x]{inputenc}
\usepackage{amssymb}
\usepackage{amsmath}
\usepackage{mathrsfs}%Serve a poter scrivere in corsivo inglese
\usepackage{bm}
\usepackage{amsthm}%Serve per avere ambienti non numerati
\usepackage{enumitem}%Serve per gli elenchi numerati
\usepackage{color}%Serve per scrivere con i colori
\usepackage{graphicx}
\usepackage{bm}%serve a scrivere in grassetto corsivo
\usepackage{ulem}%serve a barrare o sottolineare
\usepackage{upgreek}

\newcommand{\ci}[1]{\mathscr{#1}}%corsivo inglese
\newcommand{\g}[1]{\mathfrak{#1}}

\newcommand{\alfa}{\alpha}
\newcommand{\R}{\mathbf{R}}
\newcommand{\C}{\mathbf{C}}

\newcommand{\e}{\varepsilon}
\newcommand{\bra}{\left\langle}
\newcommand{\ket}{\right\rangle}
\renewcommand{\span}{\operatorname{span}}
\renewcommand{\phi}{\varphi}
\newcommand{\mi}{\mu}

\newtheorem{proposizione}{Proposition}
\newtheorem{teorema}[proposizione]{Theorem}
\newtheorem{lemma}[proposizione]{Lemma}

\DeclareMathOperator{\Sym}{Sym}
\DeclareMathOperator{\codim}{codim}

\title{Rectifiability of sets of solutions of first order systems of partial differential equations}
\author{Claudio Afeltra\footnote{Scuola Normale Superiore, Piazza dei Cavalieri 7, 56126 Pisa (Italy) - claudio.afeltra@sns.it}}
\date{}

\begin{document}

\maketitle
 
\section{Introduction}
Many mathematical problems are equivalent to impose that the differential of an unknown function $F:\Omega\to\R^m$ (with $\Omega\subseteq\R^n$) belongs to a certain subset of $\ci{L}(\R^n,\R^m)$, possibly dependent on $x\in\Omega$ and $F(x)$.

For example this is the case for holomorphic maps, isometries and conformal maps.
In some of these problems the set of solutions turns out to be finite dimensional.

In this work, given $\ci{M}\subset\Omega\times\ci{L}(\R^n,\R^m)$ we study the problem of finding conditions on $\ci{M}$ which guarantee that the set
$$\Sigma(\ci{M})=\left\{ F\in\ci{C}^1(\Omega,\R^m) \; \middle|\; \forall x\in\Omega \;\; (x,DF(x))\in\ci{M} \right\}$$
is rectifiable.
We also study a more general version of this question in which $(x,F(x),DF(x))$ is costrained to a subset of $\Omega\times\R^m\times\ci{L}(\R^n,\R^m)$, and the generalization to manifolds.

%More generally, given $\ci{M}\subset(\Omega\times\R^m)\times\ci{L}(\R^n,\R^m)$, we study the same question for the problem
%$$DF(x)\in\ci{M}_{x,F(x)} ,$$
%where $\ci{M}_{x,y}=\ci{M}\cap\pi_1(x,y)$ (with $\pi_1:(\Omega\times\R^m)\times\ci{L}(\R^n,\R^m)\to\Omega$ the projection on the first factor).

We impose that $\ci{M}$ is a manifold, and after having studied the linear case, we impose that the tangent spaces of $\ci{M}$ at the various points satisfy the conditions that we found for affine subspaces.

Let $\g{V}$ be a linear subspace of $\ci{L}(\R^n,\R^m)$.
We define the spaces of symmetric tensors
$$M_k(\g{V}) = \Big\{T\in \Sym^k(\R^n;\R^m) \;\Big|\; \forall x_1,\ldots,x_{k-1}\in\R^n$$
$$\;\text{the operator}\;x\mapsto T(x_1,\ldots,x_{k-1},x)\;\text{belongs to}\;\g{V}\Big\}$$
and
$$\g{M}(\g{V}) = \bigoplus_{k=0}^{\infty}M_k(\g{V})\subseteq \Sym(\R^n;\R^m).$$
Furthermore we define the quantities
$$\delta(\g{V})= \sup\{k \;|\; M_k\ne 0\}\in \mathbf{N}\cup\{\infty\}$$
and
$$\alpha_k(\g{V})= \dim M_k(\g{V}), \;\;\;\;    \alpha(\g{V})= \dim\g{M}(\g{V}) =\sum_{k=0}^{\infty}\alpha_k(\g{V}).$$

Obviously $\delta(\g{V})<\infty$ if and only if $\alfa(\g{V})<\infty$.

We recall that there is a natural correspondence between the space of symmetric tensors $\Sym(\R^n;\R^m)$ and the space $\R[\R^n,\R^m]$ of $\R^m$-valued polynomials on $\R^n$, given
%, on $\Sym^k(\R^n;\R^m)$, 
by $T\mapsto p(x)=T(x,\ldots,x)$. We call $\ci{P}(\g{V})$ the space of polynomials corresponding to $\g{M}(\g{V})$.

Then, in the linear case in which $\ci{M}=\Omega\times\g{V}$, we have the following theorem.

\begin{teorema}\label{CasoLineare}
 If $\ci{M}=\Omega\times\g{V}$ then $\Sigma(\ci{M})$ is a vector space which is finite dimensional if and only if $\delta(\g{V})<\infty$. In this case $\Sigma(\ci{M})$ coincides with the set of restrictions of elements of $\ci{P}(\g{V})$ to $\Omega$.
 The thesis holds also if in the definition of $\Sigma(\ci{M})$ the space $\ci{C}^1(\Omega,\R^m)$ is replaced with the space of vector-valued distributions $\ci{D}'(\Omega,\R^m)$.
\end{teorema}

Next we prove a characterization of linear subspaces $\g{V}$ of $\ci{L}(\R^n,\R^m)$ such that $\delta(\g{V})<\infty$.
There are two obvious obstructions.
If $\g{V}$ contains a rank one operator $\psi\otimes w$ with $\psi\in(\R^n)^*$ and $w\in\R^m$, then $F(x)=f(\psi(x))w\in\Sigma(\Omega\times\g{V})$ for every $f\in\ci{C}^1(\R)$.
The second obstruction is analogous, but using complex differentiable functions instead of real ones: if, up to operator equivalences, $\g{V}$ contains the space
\begin{equation}\label{DefinizioneW}
\g{W}= \span\left\{
 \left(
  \begin{array}{cc}
 	I_2 & 0 \\
 	0 & 0
  \end{array} \right),
 \left(
  \begin{array}{cc}
 	J_2 & 0 \\
 	0 & 0
  \end{array} \right)
 \right\},
\end{equation}
where $I_2 =
 \left(\begin{smallmatrix}
  1 & 0\\
  0 & 1
 \end{smallmatrix}\right)$,
 $J_2 =
 \left(\begin{smallmatrix}
  0 & -1\\
  1 & 0
 \end{smallmatrix}\right)$,
then
$$F(x)=(\g{Re}\, u(x_1+ix_2),\g{Im}\, u(x_1+ix_2),0,\ldots,0)$$
belongs to $\Sigma(\Omega\times\g{V})$ for every entire function $u$.

We will prove that these are the only two obstructions.

\begin{teorema}\label{caratterizzazione}
 $\delta(\g{V})<\infty$ if and only if $\g{V}$ does not contain neither rank one operators, nor a vector subspace of the form $P\g{W}Q$ with $P\in GL(m)$, $Q\in GL(n)$, and $\g{W}$ is the space defined in formula \eqref{DefinizioneW}.
\end{teorema}

Thanks to these results in the linear case we can prove a more general theorem.

Let $\ci{P}^*(\g{V})$ be the subspace of $\ci{P}(\g{V})$ of polynomials whose first degree term is zero.

\begin{teorema}\label{Teorema}
 Let $\ci{M}$ be a submanifold of $\Omega\times\ci{L}(\R^n,\R^m)$ of class $\ci{C}^1$ such that the projection $\pi_1$ on $\Omega$ restricted to $\ci{M}$ is a submersion.
 
 Calling $\ci{M}_x=\ci{M}\cap\pi_1^{-1}(x)$ and $\g{V}_{x,A}=T_A(\ci{M}_x)-A$, let us suppose that $\alpha(\g{V}_{x,A}) =k$ for every $x\in\Omega$ and $A\in\ci{M}_x$, and that $\Omega$ is connected.
 
 Then $\Sigma(\ci{M})$ is $k$-rectifiable.
 
 Furthermore if $F,G\in\Sigma(\ci{M})$ coincide on a non empty open subset of $\Omega$, they coincide on $\Omega$.
 
 Finally the functions $\alpha_j$ and $\delta$ are constant, and calling $\delta(\g{V}_{x,A})=d$, if $\widetilde{\Sigma}\subseteq\Sigma(\ci{M})$ is a $\ci{C}^1$ submanifold of $\ci{C}^1(\Omega,\R^m)$ such that for some $x\in\Omega$ every function in $\widetilde{\Sigma}$ is of class $\ci{C}^{d}$ in a neighborhood of $x$ and for every $A\in\ci{M}_x$ and $p\in\ci{P}^*(\g{V}_{x,A})$ there exists $F\in\widetilde{\Sigma}$ with Taylor expansion corresponding to $A$ and $p$ then $\widetilde{\Sigma}=\Sigma(\ci{M})$.
\end{teorema}

Finally we treat the more general case in which the set to which $DF$ is costrained depends also on $F$. This allows to state and prove the theorem in the more general setting of manifolds.

Let $N,M$ be two manifolds with $\dim N=n$, $\dim M=m$, and consider the bundle $\pi_1^*T^*N\otimes\pi_2^*TM$ on $N\times M$. Let $p:\pi_1^*T^*N\otimes\pi_2^*TM\to N\times M$ be the projection.
Let $\ci{M}$ be a $\ci{C}^1$ submanifold of $\pi_1^*T^*N\otimes\pi_2^*TM$ such that $\pi_1\circ p|_{\ci{M}}$ is a submersion, and let
$$\Sigma(\ci{M}) = \left\{ F\in\ci{C}^1(N,M) \;\middle|\;\forall x\in N\right.$$
$$\left.d_xF\in(\pi_1^*T^*N\otimes\pi_2^*TM)_{x,F(x)}\text{ belongs to }\ci{M}\right\}.$$
Let
$$\ci{M}_{x,y}=\ci{M}\cap p^{-1}(x,y)\subset\ci{L}(T_xN,T_yM),$$
and given $A\in\ci{M}_{x,y}$, let
$$\mathbf{V}_{x,y,A}\subset (T_x^*N\otimes T_yM)\oplus T_yM\simeq\ci{L}(T_xN,T_yM)\oplus T_yM$$
be the subspace corresponding to the linearization of the problem in $A$, and $\g{V}_{x,y,A}\subset\ci{L}(T_xN,T_yM)$ its projection on the first factor.

Let
$$\mathbf{P}_{x,y,A}= \big\{u:T_xN\to T_yM \;\big|\; u(0)=0, Du(0)=A, $$
$$\forall\xi\in T_xN \;(Du(\xi),u(\xi))\in\mathbf{V}_{x,y,A}\big\}.$$

In the Euclidean case, $\pi_1^*T^*N\otimes\pi_2^*TM = (\Omega\times\R^m)\times\ci{L}(\R^n,\R^m)$, $p$ is the projection on $\Omega\times\R^m$ and
$$\Sigma(\ci{M})=\left\{ F\in\ci{C}^1(\Omega,\R^m) \; \middle|\; \forall x\in\Omega \;\; (x,F(x),DF(x))\in\ci{M} \right\}.$$

%let $\ci{P}^{**}(\g{V})$ be the subspace of $\ci{P}(\g{V})$ with zero constant and first order term.

Then the following result holds.

\begin{teorema}\label{Teorema2}
 In the above setting, suppose that $\alpha(\g{V}_{x,y,A}) =k$ for every $x\in N$, $y\in M$ and $A\in\ci{M}_x$, and that $N$ is connected.
 
 Then $\Sigma(\ci{M})$ is $k$-rectifiable.
 
 Furthermore if $F,G\in\Sigma(\ci{M})$ coincide on a non empty open subset of $N$, they coincide on $N$.
 
 Finally $\dim\mathbf{P}_{x,y,A}=k-\dim M-\dim\g{V}_{x,y,A}$, 
 there exists $d$ such that the elements of $\mathbf{P}_{x,y,A}$ have different Taylor expansions of order $d$ at zero, and if 
 %if $\delta(\g{V}_{x,y,A})=d$,
 if $\widetilde{\Sigma}\subseteq\Sigma(\ci{M})$ is a manifold of class $\ci{C}^1$ is such that for some $x\in N$
 every function in $\widetilde{\Sigma}$ is of class $\ci{C}^{d}$ in a neighborhood of $x$
 and for every $y\in M$, $A\in\ci{M}_{x,y}$
 %and $p\in\ci{P}^{**}(\g{V}_{x,y,A})$
 and $u\in\mathbf{P}_{x,y,A}$ there exists $F\in\Sigma$ such that $F(x)=y$ with Taylor expansion of order $d$ corresponding to $u$,
 %the set
 %$$\widetilde{\Sigma}_{x,y,A}=\left\{ F\in\widetilde{\Sigma} \;\middle|\; F(x)=y, DF(x)=A \right\}$$
 %is a $\ci{C}^1$ manifold of dimension $\alfa(\g{V}_{x,y,A})-\dim M-\dim\ci{M}_x$, 
 then $\widetilde{\Sigma}=\Sigma(\ci{M})$.
\end{teorema}

In the last section, we show some examples of application of our results and find alternative proofs to some known facts.
%All of them are known facts, but we believe that the alternative proofs we show are interesting and 

\section{The linear case}
\begin{proof}[Proof of Theorem \ref{caratterizzazione}]
 Let us suppose that $\delta(\g{V})=\infty$.
 Let us define
 $$(\Lambda_xT)(x_1,\ldots,x_{k-1})= T(x_1,\ldots,x_{k-1},x).$$
 
 Let us consider the complexifications of $M_k(\g{V})$, $M_k(\g{V})_{\C}$.
 
 Then $M_k(\g{V})_{\C}$ is a family of linear spaces such that
 \begin{itemize}
  \item $M_k(\g{V})_{\C}\subseteq \Sym(\C^n;\C^m)$;
  \item for every $T\in M_k(\g{V})_{\C}$ and $x\in\C^n$, $\Lambda_x(T)\in M_{k-1}(\g{V})_{\C}$;
  \item if $\Lambda_xT\in M_{k-1}(\g{V})_{\C}$ for every $x\in\C^n$ then $T\in M_k(\g{V})_{\C}$;
  \item $M_k(\g{V})_{\C}\ne 0$ for every $k$.
 \end{itemize}
 A standard application of Zorn's lemma permits to prove the existence of minimal family of spaces $\widetilde{M}_k\subseteq M_k(\g{V})_{\C}$ such that:
 \begin{enumerate}[label = (\roman*)]
  \item\label{CondizioneI} for every $T\in \widetilde{M}_k$ and $x\in\C^n$, $\Lambda_x(T)\in\widetilde{M}_{k-1}$;
  \item\label{CondizioneII} $\widetilde{M}_k\ne 0$ for every $k$.
 \end{enumerate}
 In fact the set of family of spaces verifying the above conditions is partially ordered by inclusion.
 If $\{M_k^{\alfa}\}_{\alfa\in\ci{A}}$ is a chain indexed by a totally ordered set $\ci{A}$ then $M'_k=\bigcap_{\alfa\in\ci{A}}M_k^{\alfa}$ verifies condition \ref{CondizioneI} trivially and condition \ref{CondizioneII} because the intersection of a chain of non-zero vector subspaces of a finite dimensional vector subspaces is non-zero. Thus Zorn's lemma applies.
 
 For every fixed $x\in\C^n$ the family $\left(\Lambda_x\widetilde{M}_{k+1}\right)_{k\in\mathbf{N}}$ verifies condition \ref{CondizioneI} and $\Lambda_x\widetilde{M}_{k+1}\subseteq \widetilde{M}_k$, thus, by minimality, either it coincides with $\widetilde{M}_k$ or it does not verify condition \ref{CondizioneI}, that is $\Lambda_x\widetilde{M}_{k+1}=0$ eventually.
 The set
 $$V=\left\{ x\in\C^n \;\middle|\; \Lambda_x\widetilde{M}_{k+1}=0 \;\text{eventually}\right\}$$
 is obviously a subspace of $\C^n$, and it is proper because if $\Lambda_xT=0$ for every $x\in\C^n$ then $T=0$, and so if $V=\C^n$ then condition \ref{CondizioneII} would be contradicted.
 For $x\in\C^n\setminus V$ it holds that $\Lambda_x\widetilde{M}_{k+1}=\widetilde{M}_k$.
 
 Analogously
 %if $N_{x,k}= \ker{\Lambda_x|_{M_k}}$ then 
 the family $\left(\ker{\Lambda_x|_{\widetilde{M}_k}}\right)_{k\in\mathbf{N}}$ verifies condition \ref{CondizioneI} and is such that $\ker{\Lambda_x|_{\widetilde{M}_k}}\subseteq\widetilde{M}_k$, so it is either eventually zero ot it coincides with $\widetilde{M}_k$.
 The set
 $$W=\left\{ x\in\C^n \;\middle|\; \forall k \;\; \ker{\Lambda_x|_{\widetilde{M}_k}}=\widetilde{M}_k \right\}$$
 %and is eventually zero for $x$ in a proper subspace of $\C^n$.
 is a subspace of $\C^n$, and it is proper because $\widetilde{M}_k$ satisfies condition \ref{CondizioneII}.
 For $x\in\C^n\setminus W$, $\ker{\Lambda_x|_{\widetilde{M}_k}}=0$ eventually.
 
 In particular for $x\in(\C^n\setminus V)\cap(\C^n\setminus W)=\C^n\setminus(V\cup W)$ (which is obviously non empty)
 %outside the union of two proper subspaces
 and $k$ big enough, $\Lambda_x:\widetilde{M}_k\to \widetilde{M}_{k-1}$ is an isomorphism.
 So $\dim\widetilde{M}_k=\dim \widetilde{M}_{k-1}$ and $\Lambda_x:\widetilde{M}_k\to \widetilde{M}_{k-1}$ is injective if and only if it is surjective, hence $\C^n\setminus V=\C^n\setminus W$, that is $V=W$.
 
 Let us suppose by contradiction that $\codim V>1$.
 Let $x,y$ in some complementary of $V$. Then $\Lambda_{x+\xi y}:\widetilde{M}_k\to \widetilde{M}_{k-1}$ is an isomorphism for every $\xi\in\C$. But $\det(\Lambda_{x+\xi y})=\det(\Lambda_x+\xi \Lambda_y)$ (determinant with respect to two bases of $\widetilde{M}_k$ and $\widetilde{M}_{k-1}$) is a non constant polynomial, and this is absurd. So $\codim V=1$. 
 
 %If $\Lambda_x$ and $\Lambda_y$ are linearly independent then $\Lambda_x+\xi \Lambda_y$ is an isomorphism for every $\xi$. But $\det(\Lambda_x+\xi \Lambda_y)$ (determinant with respect to two bases of $\widetilde{M}_k$ and $\widetilde{M}_{k-1}$) is a non constant polynomial, and this is absurd. So the operators $L_x$ form a space of dimension one.
 
 Let $\psi\in(\C^n)^*$ such that $V=\ker\psi$. Then if $(\psi_1=\psi,\psi_2,\ldots,\psi_n)$ is a basis of $(\C^n)^*$ which is the dual basis of $(w_1,\ldots,w_n)$, and $k$ is big enough such that $\Lambda_x\widetilde{M}_k=0$ for $x\in V$, then testing this condition for $x=w_2,\ldots,w_n$, it is easily proved that $\widetilde{M}_k$ is generated by $T=\psi\otimes\ldots\otimes\psi\otimes v$ for some $v\in\R^m$.
 
 Therefore $(\Lambda_{w_1})^{k-1}T=\psi\otimes v\in\widetilde{M}_1\subseteq M_1(\g{V})_{\C}=\g{V}\otimes\C$.
 Now the thesis follows by taking the real and imaginary part of $\psi\otimes v$, the two cases stated in the theorem arising respectively whether the real and imaginary part thereof be linearly dependent or not.
\end{proof}

\begin{proof}[Proof of Theorem \ref{CasoLineare}]
 Let $F\in\Sigma(\ci{M})$. By convolution one can suppose in advance that $F$ is smooth.
 Then the $k$-th differential of $F$, $D^k F$, is a function from $\Omega$ to $M_k(\g{V})$. From this the thesis follows easily.
\end{proof}

\section{The general case}
Let $\g{G}_k(V)$ be the $k$-Grassmannian of the vector space $V$.

\begin{lemma}
 For every $\ell\in\mathbf{N}$ the function $\dim M_{\ell}:\g{G}_k(\ci{L}(\R^n,\R^m))\to \mathbf{N}\cup\{\infty\}$ is upper semicontinuous, and thus so is their sum $\alpha$. Given $d\in\mathbf{N}$, the function $\ci{P}:\alpha^{-1}(d)\to \g{G}_d(\R[\R^n;\R^m])$ is continuous. An analogous statement holds for $\ci{P}^*$.
\end{lemma}

\begin{proof}
 Let $\g{V}_j$ a sequence of $k$-dimensional subspaces of $\ci{L}(\R^n,\R^m)$ with
 \\
 $\dim M_{\ell}(\g{V}_j)=\mi$ converging to $\g{V}$. $M_{\ell}(\g{V}_j)$ is a sequence in $\g{G}_{\mi}(\Sym^k(\R^n;\R^m))$, and so by the compactness of the Grassmannians it has a subsequence convergent to some $W\in\g{G}_{\mi}(\Sym^k(\R^n;\R^m))$. By continuity for every $T\in W$ and $x_1,\ldots,x_{\ell-1}\in\R^n$ the operator $x\mapsto T(x_1,\ldots,x_{\ell-1},x)$ belongs to $\g{V}$, and so $W\subseteq M_{\ell}(\g{V})$, which implies that $\dim M_{\ell}(\g{V})\ge\mi$.
 
 If $\alpha(\g{V}_j)= \alpha(\g{V})=d$ then by semicontinuity $\dim M_{\ell}(\g{V}_j) \le \dim M_{\ell}(\g{V})$ eventually in $j$ for every $\ell$, and since $\sum\dim M_{\ell}(\g{V}_j)=\alpha(\g{V}_j)= \alpha(\g{V})=\dim M_{\ell}(\g{V})$, then $\dim M_{\ell}(\g{V}_j) = \dim M_{\ell}(\g{V})$ eventually. Thus, repeating the above reasoning, by dimensional reasons $W=M_{\ell}(\g{V})$. Thus every subsequence of $M_{\ell}(\g{V}_j)$ has a subsequence convergent to $M_{\ell}(\g{V})$, and therefore the second statement of the lemma follows.
\end{proof}

\begin{proof}[Proof of Theorem \ref{Teorema}]
 Let $\phi:\Omega\times\ci{L}(\R^n,\R^m)\to\R^{nm-k}$ be a defining function for $\ci{M}$, and let 
 $\Phi(u)(x)=\phi(x,Du(x))$,
 $$\Phi:\ci{C}^1(\Omega,\R^m)\to\ci{C}(\Omega,\R^{nm-k})$$
 and let $\Phi_U:\ci{C}^1(U,\R^m)\to\ci{C}(U,\R^{nm-k})$ for $U\subseteq\Omega$ be defined analogously.

 Then the set to be studied is $\Phi^{-1}(0)$.
 
 $\Phi_U$ is of class $\ci{C}^1$ with
 $$(D\Phi_U(u)[v])(x) = D_2\phi(x,Du(x))[Dv(x)].$$
 Let
 $$\ci{C}^1_{x_0}(U,\R^m) = \left\{u\in\ci{C}^1(U,\R^m) \;\middle|\; Du(x_0)=0 \right\}.$$
 We know that, given $x_0\in\Omega$, if
 $$A_{u,x_0}= D_2\phi(x_0,Du(x_0))\in\ci{L}(\ci{L}(\R^n,\R^m),\R^k)$$
 and for $U$ in a neighborhood of $x_0$ and $v\in\ci{C}^1_{x_0}(U)$
 $$(L_{u,x_0,U}[v])(x) = A_{u,x_0}[Dv(x)],$$
 then $L_{u,x_0,U}:\ci{C}^1_{x_0}(U,\R^m)\to\ci{C}(U,\R^k)$ has a finite dimensional kernel, equal to $\ci{P}^*(\ker A_{u,x_0})=\ci{P}^*(\g{V}_{x_0,Du(x_0)})$
\footnote{Technically it consists of the restrictions of $\ci{P}^*(\g{V}_{x_0,Du(x_0)})$ to $U$, but obviously we can identify it with $\ci{P}^*(\g{V}_{x_0,Du(x_0)})$ independently from $U$}.

 Let $\pi:\ci{C}^1_{x_0}(U,\R^m)\to\ker L_{u,x_0,U}$ be a projection. Then $L_{u,x_0,U}$ is invertible on $\ker\pi$. By continuity if $U$ is small enough then $D\Phi_U(u)$ is invertible on $\ker\pi$, and so also $\Phi_U$ is.

 Suppose that $F,G\in\Phi^{-1}(0)$ coincide on an open set, but not on the whole $\Omega$. Then there exists $x_0$ belonging to the boundary of the inner part of $\{F=G\}$.
 Since $x_0$ belongs to the clousure of $\{F=G\}$ there exist $y$ arbitrarily close to $x_0$ such that $B_{\e}(y)\subset \{F=G\}\cap B_r(x_0)$ for some $\e$.
 %Let us consider a projection $\pi:\ci{C}^1_{y}(B_{\e}(y),\R^m)\to\ker L_{F,y,B_{\e}(y)}=\ci{P}^*(\g{V}_{y,DF(y)})$.

 Let us consider a projection $\pi:\ci{C}^1(B_{\e}(y),\R^m)\to\ci{P}(\g{V}_{y,DF(y)})$.
 Obviously $\pi((F-G)|_{B_{\e}(y)})=0$. By composing the inclusion
 $$\ci{C}^1_{x_0}(B_r(x_0),\R^m)\to\ci{C}^1(B_r(x_0),\R^m),$$
 the restriction to $B_{\e}(y)$, $\pi$ and the natural projection
 $$\ci{P}(\g{V}_{y,DF(y)})\to\ci{P}^*(\g{V}_{y,DF(y)}),$$
 we get a projection $\pi':\ci{C}^1_{x_0}(B_r(x_0),\R^m)\to\ci{P}^*(\g{V}_{y,DF(y)})$. By continuity its restriction to $\ci{P}^*(\g{V}_{x_0,DF(x_0)})$ is an isomorphism if $|y-x_0|$ is small enough, and so by composing with an automorphism of $\ci{C}^1_{x_0}(B_r(x_0),\R^m)$ which on $\ci{P}^*(\g{V}_{x_0,DF(x_0)})$ coincides with the inverse thereof, we get a projection
 $$\pi'':\ci{C}^1_{x_0}(B_r(x_0),\R^m)\to\ci{P}^*(\g{V}_{x_0,DF(x_0)})$$
 such that $\pi''(F-G)=0$.
 Since, by the above reasoning, $\Phi_{B_r(x_0)}$ is invertible on $\ker\pi''$, $F=G$ in $B_r(x_0)$, and this is a contradiction.

 So we can restrict to arbitrarily small neighborhoods of some point $x_0$.
 Let us consider $F\in\Sigma(\ci{M})$. By the above reasoning there exist a cone in $\ci{C}^1_{x_0}(U)$ around $\ci{P}^*(\g{V}_{x_0,DF(x_0)})$ containing $F-G$ for any element $G$ of $\Sigma(\ci{M})$ in a neighborhood of $F$ such that $DG(x_0)=DF(x_0)$.
 
 Indeed
 $$0=\Phi(G)-\Phi(F)= D\Phi(F)[G-F] + o(G-F) = $$
 $$=D\Phi(F)[\pi(G-F)] - D\Phi(F)[G-F-\pi(G-F)]+ o(G-F) = $$
 $$= - D\Phi(F)[G-F-\pi(G-F)]+ o(G-F)$$
 and $G-F-\pi(G-F)\in\ker\pi$, on which $D\Phi(F)$ is invertible.
 
 By continuity there exists a cone in $\ci{C}^1(U)$ containing any element of $\Sigma(\ci{M})$ in a neighborhood of $F$. By Theorem 6.2 in \cite{AK} (which can be applied to $\ci{C}^1(\Omega,\R^m)$ because it is isometrically embedded in $W^{1,\infty}(\Omega,\R^m)$ which is the dual of a separable Banach space), $\Sigma(\ci{M})$ is $k$-rectifiable.

 As for the last part of the theorem, defining
 $$\widetilde{\Sigma}_F=\left\{G\in\widetilde{\Sigma}\;\middle|\;DG(x)=DF(x)\right\},$$
 in the above reasoning one can choose a projection $\pi$ mapping $\widetilde{\Sigma}_F-F$ bijectively to $\ci{P}^*(\g{V}_{x,DF(x)})$, and thus the latter rigidity stamement allows to conclude.
\end{proof}

\begin{proof}[Proof of Theorem \ref{Teorema2}]
 The theorem is proved analogously to Theorem \ref{Teorema}. The only substantial modification is in regard to $\mathbf{P}_{x,y,A}$, whose study is reduced to $\ci{P}(\g{V}_{x,y,A})$ by known methods of linear systems of partial differential equations (see \cite{O}).
\end{proof}

\section{Applications}
%In this section we show some applications of the above results. They are all already known facts, but we believe that this shows the usefulness of the work done.

\subsection{Conformal maps}
Theorem \ref{Teorema} allows to prove the Liouville theorem stating that conformal maps between open sets of $\R^n$ with $n\ge 3$ are compositions of isometries and spherical inversions.
In fact in this case $\ci{M}=\Omega\times CO(n)$, where $CO(n)$ is the conformal group
$$CO(n)= \left\{ A\in GL(n) \;\middle|\; AA^T= (\det A)^2I\right\}.$$
Since $CO(n)$ is a Lie group, $\alfa(\g{V}_{x,A})= \alfa(\g{V}_{x,I})$. It is easily verified that
$$\g{V}_{x,I} = \left\{ \lambda I +A \;\middle|\; \lambda\in\R, A+A^T=0 \right\}.$$
The hypothesis of Theorem \ref{caratterizzazione} are verified if $n\ge 3$ (and are not when $n=2$, but in that case conformal maps are holomorphic or anti-holomorphic maps with nowhere zero derivative, and thus form an infinite-dimensional set).

It is obvious that $M_0(\g{V}_{x,I})=\R^n$ and $M_1(\g{V}_{x,I})=\g{V}_{x,I}$.

If $B\in M_2(\g{V}_{x,I})$ then there exist $\alfa\in(\R^n)^*$ and $A:\R^n\to\ci{L}(\R^n)$ such that $A(x)$ is skew-symmetric and $B(x,y)=\alfa(x)y+A(x)y$. Therefore
$$\bra B(x,y),z \ket = \alfa(x)\bra y,z\ket+\bra A(x)y,z\ket = \alfa(x)\bra y,z\ket-\bra y,A(x)z\ket =$$
$$= 2\alfa(x)\bra y,z\ket-\bra B(x,z),y \ket =2\alfa(x)\bra y,z\ket-\bra B(z,x),y \ket =$$
$$=2\alfa(x)\bra y,z\ket-\alfa(z)\bra x,y\ket - \bra A(z)x,y \ket$$
and thus
$$\bra B(x,y),z \ket =\bra B(y,x),z \ket =2\alfa(y)\bra x,z\ket-\alfa(z)\bra x,y\ket - \bra A(z)y,x \ket=$$
$$= 2\alfa(y)\bra x,z\ket-\alfa(z)\bra x,y\ket + \bra A(z)x,y \ket$$
and summing we get that
\begin{equation}\label{BilineariLiouville}
 \bra B(x,y),z\ket= \alfa(x)\bra y,z\ket +\alfa(y)\bra x,z\ket -\alfa(z)\bra x,y\ket.
\end{equation}
Thus $M_2(\g{V}_{x,I})$ consists of elements of $\Sym^2(\R^n;\R^n)$ of the above form \eqref{BilineariLiouville}.

If $T\in M_3(\g{V}_{x,I})$ then there exists a bilinear form $b$ such that
$$\bra T(x,y,w),z\ket = b(x,w)\bra y,z\ket +b(y,w)\bra x,z\ket -b(z,w)\bra x,y\ket;$$
therefore
$$\bra T(y,z,w),x\ket = b(y,w)\bra z,x\ket +b(z,w)\bra y,x\ket -b(x,w)\bra y,z\ket,$$
$$\bra T(z,x,w),y\ket = b(z,w)\bra x,y\ket +b(x,w)\bra z,y\ket -b(y,w)\bra z,x\ket$$
and the last three equations imply that $b=0$, and thus $M_3(\g{V}_{x,I})=0$.

Hence we found $\g{M}(\g{V}_{x,I})$, which is equivalent to find $\ci{P}(\g{V}_{x,I})$. Now to prove Liouville theorem it is sufficient to prove that the set $\widetilde{\Sigma}$ given by compositions of isometries and spherical inversions verifies the hypotheses of Theorem \ref{Teorema}.

The above computations work also for conformal maps between general $n$-dimensional Riemannian manifolds with $n\ge 3$, and allows to prove that the set of conformal maps between two manifolds is $\frac{(n+1)(n+2)}{2}$-rectifiable, and in a similar manner one can prove that the set of isometries is $\frac{n(n+1)}{2}$-rectifiable. Obviously both facts are much weaker than the known results that both sets are manifolds, the former of dimension at most $\frac{(n+1)(n+2)}{2}$ and the latter of dimension at most $\frac{n(n+1)}{2}$ (see \cite{K}).

\subsection{Quaternion differentiable functions}
Let $\mathbf{H}$ be the set of quaternions. A function $f:\mathbf{H}\to\mathbf{H}$ is said left quaternion differentiable if the limit
$$\frac{df}{dq} = \lim_{h\to 0}h^{-1}(f(q+h)-f(q))$$
exists. This is equivalent to $f$ being real differentiable, and the differential being the quaternion multiplication to the right by $\frac{df}{dq}$.

It is known that left quaternion differentiable functions are actually affine. We can prove this fact, under the additional hypothesis that $f$ is $\ci{C}^1$, by applying our theory: in this case $\ci{M}=\Omega\times\g{V}$ where $\g{V}$ is the subspace of $\ci{L}(\mathbf{H})$ given by right multiplications by some quaternion. Since the non zero elements of $\g{V}$ have rank four, $\g{V}$ satisfies the hypothesis of Theorem \ref{caratterizzazione}.

If $B\in M_2(\g{V})$ then there exists an operator $A\in\ci{L}(\mathbf{H})$ such that $B(x,y)=y\cdot Ax$. Taking $y=1$ and using symmetry,
$$Ax=1\cdot Ax=B(x,1)=B(1,x)=x\cdot A1.$$
Defining $a=A1$, $B(x,y)=y\cdot x\cdot a$. By symmetry $y\cdot x\cdot a= B(x,y)=B(y,x)=x\cdot y\cdot a$ for every $x,y\in\mathbf{H}$, and this is possible only if $a=0$. Thus $M_2(\g{V})=0$ and $\Sigma(\Omega\times\g{V})$ consists only of affine functions.

\subsection{Local CR isomorphisms}
We refer to \cite{DT} for basic concepts in CR geometry.

%Let $\Omega$ be an open subset of $S^{2n+1}$ with its standard CR structure.
Let $N$ and $M$ be two $2n+1$-dimensional CR manifolds.
We want to study local CR isomorphisms
%$\Omega\to S^{2n+1}$
$F:N\to M$
with Theorem \ref{Teorema2}. If, given
%$x\in\Omega$
$x\in N$ and
%$y\in S^{2n+1}$
$y\in M$
we define
%$$\ci{M}_{x,y}= \left\{ A:T_x\Omega\to T_yS^{2n+1} \;\middle|\; A(T^{(1,0)}_x\Omega)=T^{(1,0)}_yS^{2n+1} \right\},$$
$$\ci{M}_{x,y}= \left\{ A:T_xN\to T_yM \;\middle|\; A(T^{(1,0)}_xN)=T^{(1,0)}_yM \right\},$$
then the seeked isomorphisms are $\Sigma(\ci{M})$, but the manifold $\ci{M}$ does not satisfy the hypotheses of Theorem \ref{Teorema2}.

If we impose the condition that the differential of a CR isomorphism is conformal with respect to the Levi form, and that (fixing two contact forms $\theta_N$ and $\theta_M$ on the two CR structures) if $F^*L_{\theta_M}=uL_{\theta_N}$ then $\theta_M(dF(T_N)) = u$ (where $T_N$ is the Reeb vector field with respect to $\theta_N$), we can describe local CR isomorphisms as $\Sigma(\ci{M}')$ where
$$\ci{M}'_{x,y}= \Big\{ A:T_xN\to T_yM \;\Big|\; A(T^{(1,0)}_xN)=T^{(1,0)}_yM,$$
$$\exists\lambda>0 \;\text{such that} \; A^*((L_{\theta_M})_y)=\lambda (L_{\theta_N})_x \;\text{and}\; (\theta_M)_y(A(T_N)_x)=\lambda \Big\},$$
but yet again $\ci{M}'$ does not verify the hypotheses of Theorem \ref{Teorema2}. The problem is that in the definition of $\ci{M}'$ no condition is imposed on the component of $dF(T_N)$ along $H(M)$, and because of this $\g{V}'_{x,y}$ contains rank one operators.
To solve the problem we note that by the formula for the conformal transformation of the Reeb field the component of $(dF(T_N))_y$ along $H_y(M)$ is determined by the restriction of $dF$ to $T^{(1,0)}_xN$. This allows to build a third manifold $\ci{M}''$ such that local CR isomorphisms are described as $\Sigma(\ci{M}'')$ and it verifies the hypotheses of Theorem \ref{Teorema2}. Then computations very similar to the ones for the conformal maps performed above allow to prove that the set of local CR isomorphisms is rectifiable, and the well known fact that in the case of $S^{2n+1}$ the group of restrictions of biholomorphisms of the ball of $\C^{n+1}$, which is naturally identified with $PU(n+1,1)$.

This example is interesting because it shows that in some cases the curvature of $\ci{M}$ can ensure that the thesis of Theorem \ref{Teorema2} holds even if $\ci{M}$ does not verify the hypotheses thereof.

\end{document}